\documentclass{article}%
\usepackage{amsfonts}%
\usepackage{amsmath}%
\usepackage{amssymb}%
\usepackage{graphicx}
\providecommand{\U}[1]{\protect\rule{.1in}{.1in}}
\newtheorem{theorem}{Theorem}

\newtheorem{corollary}[theorem]{Corollary}

\newtheorem{lemma}[theorem]{Lemma}

\newtheorem{proposition}[theorem]{Proposition}

\ifx\pdfoutput\relax\let\pdfoutput=\undefined\fi
\newcount\msipdfoutput
\ifx\pdfoutput\undefined\else
\ifcase\pdfoutput\else
\ifx\paperwidth\undefined\else
\ifdim\paperheight=0pt\relax\else\pdfpageheight\paperheight\fi
\ifdim\paperwidth=0pt\relax\else\pdfpagewidth\paperwidth\fi
\fi\fi\fi
\begin{document}

\title{ $H^{o}$-type Riemannian metrics on the space of planar curves \thanks{This
work was supported by NIH Grant I-R01-NS34189-08}}
\author{\bigskip
\and \bigskip Jayant Shah\\Mathematics Department, \ Northeastern University, Boston, MA\\email: shah@neu.edu}
\maketitle

\begin{abstract}
Michor and Mumford have shown that the distances between planar curves in the
simplest metric (not involving derivatives) are identically zero. We consider
two conformally equivalent metrics for which the distances between curves are
nontrivial. We show that in the case of the simpler of the two metrics, the
only minimal \ geodesics are those corresponding to curve evolution in which
the points of the curve move with the same normal speed. An equation for the
geodesics and a formula for the sectional curvature are derived; a necessary
and sufficient condition for the sectional curvature to be bounded is given.

\end{abstract}

\section{Introduction}

The purpose of this paper is to study the most basic properties of some of the
simplest Riemannian metrics suggested by applications to Computer Vision. The
problem is to understand and quantify similarities and differences between
object shapes and their individual variations. At a fundamental level, the
problem is to construct appropriate metrics on a space of closed surfaces in
$\mathbb{R}^{3}$. A simpler version of the problem is the construction of
Riemannian metrics on a space of closed planar curves. The choice of a metric
depends on the type of similarity that is being considered. In their seminal
paper [2], Michor and Mumford analyze two Riemannian metrics on a space of
closed planar curves. Surprisingly, the Riemannian distance between any two
curves in the simpler of the two metrics, an $H^{o}$-metric, turns out to be
zero. To remedy this, they add a curvature term to the metric. Alternative
construction of Riemannian metrics on the space of curves of constant length
and parametrized by either the orientation of the tangent vector or the
curvature is given by Klasssen et al in [3]. A method for deforming one curve
onto another by minimizing an approximate Hausdorff distance between them is
described in [1]. Below, we analyse two conformal variants of the $H^{o}%
$-metric of Michor and Mumford. In \S 3, we derive upper and lower bounds for
distances between curves and show that these metrics behave like $L^{1}$
metrics. A key issue is the existence of minimal geodesics. We show that for
the simpler of the two conformal metrics, the minimal geodesics correspond
exactly to those curve deformations in which the points of the curve move with
the same normal speed. In the case of the second metric, no geodesic is
minimal if the length of the curve is less than a certain threshold; the
question of minimality when the length of the deforming curve is equal or
greater than threshold is still open. We provide a partial answer\ in the form
of a necessary and sufficient condition for the boundedness of the sectional
curvature in \S 5. An equation for the geodesics is derived in \S 4.

This paper began with an analysis of the simpler of the two conformal metrics
considered here. As it was being written, the author became aware of the work
of A. Yezzi and A. Mennucci [4,5] in which they have proposed a more general
formulation. Since our analysis applies to their formulation as well, we have
included their formulation in the analysis below.

\section{The Framework}

The basic space considered by Michor and Mumford is the orbit space
\[
B_{e}(S^{1},\mathbb{R}^{2})=Emb(S^{1},\mathbb{R}^{2})/Diff(S^{1})
\]
of the space of all $C^{\infty}$ embeddings of $S^{1}$ in the plane, under the
action by composition from the right by diffeomorphisms of the unit circle. It
is contained in the bigger space of immersions modulo diffeomorphisms:
\[
B_{i}(S^{1},\mathbb{R}^{2})=Imm(S^{1},\mathbb{R}^{2})/Diff(S^{1})
\]
Let $\pi:Imm(S^{1},\mathbb{R}^{2})\rightarrow B_{i}(S^{1},\mathbb{R}^{2})$ be
the canonical projection. The simpler of the two metrics considered in [2] is
an $H^{o}$-metric defined on $Imm(S^{1},\mathbb{R}^{2})$:
\begin{equation}
G_{c}^{o}(h,k)=\int_{S^{1}}\left(  h\cdot k\right)  |c_{\theta}|d\theta
\end{equation}
where $c:S^{1}\rightarrow\mathbb{R}^{2}$ is an immersion, defining a point in
$Imm(S^{1},\mathbb{R}^{2})$, $h,k$ $\in C^{\infty}(S^{1},\mathbb{R}^{2})$ are
the vector fields along the image curve, defining two tangent vectors in
$Imm(S^{1},\mathbb{R}^{2})$ at $c$, and $c_{\theta}=dc/d\theta$. $(h\cdot k)$
is the usual dot product in $\mathbb{R}^{2}$. Sometimes for the sake of
clarity, we will use the notation $a\cdot b$ even when $a,b$ are scalars. Let
$n_{c}$ denote the unit normal field along $c$. If we identify $\mathbb{R}
^{2}$ with the complex plane $\mathbb{C}$, then, $n_{c}=ic_{\theta}/\left|
c_{\theta}\right|  $. For any $C_{o},C_{1}\in$ $B_{i}$, consider all liftings
$c_{o},c_{1}$ to $Imm(S^{1},\mathbb{R}^{2})$ and all smooth paths $t$
$\mapsto(\theta\mapsto c(t,\theta))$, $0\leq t\leq1$, in $Imm(S^{1}%
,\mathbb{R}^{2})$ with $c(0,\cdot)=c_{o}$ and $c(1,\cdot)=c_{1}$. Let $c_{t}$
denote $\partial c/\partial t$ and $c_{t}^{\bot}=\left(  c_{t}\cdot
n_{c}\right)  n_{c}$. The arc-length of such a path $c$ is given by
\[
\int_{0}^{1}\sqrt{G_{c}^{o}\left(  c_{t},c_{t}\right)  }dt
\]
Michor and Mumford show that for any two curves in $B_{i}(S^{1},\mathbb{R}%
^{2})$ ,
\[
dist_{G^{o}}(C_{1},C_{2})=_{def}\inf_{c}\int_{0}^{1}\sqrt{G_{c}^{o}\left(
c_{t}^{\bot},c_{t}^{\bot}\right)  }dt=0
\]
and strengthen $G^{o}$ by defining
\[
G_{c}^{A}(h,k)=\int_{S^{1}}\left(  1+A\kappa_{c}^{2}\right)  \left(  h\cdot
k\right)  |c_{\theta}|d\theta
\]
where $\kappa_{c}$ is the curvature, defined by the equation $\left(  \frac{
c_{\theta}}{\left|  c_{\theta}\right|  }\right)  _{\theta}=i\kappa
_{c}c_{\theta}=\kappa_{c}\left|  c_{\theta}\right|  n_{c}$.

An alternative is to consider conformal transformations of $G^{o}$. They have
the form
\[
G_{c}^{\Phi}(h,k)=\Phi(c)\int_{S^{1}}\left(  h\cdot k\right)  |c_{\theta
}|d\theta
\]
where $\Phi$ is a $Diff(S^{1})$-invariant function on $Imm(S^{1},\mathbb{R}
^{2})$. For example,
\[
\Phi(c)=\int_{S^{1}}|c_{\theta}|d\theta\text{\qquad or \qquad}\Phi
(c)=\int_{S^{1}}\left(  1+A\kappa_{c}^{2}\right)  |c_{\theta}|d\theta
\]
Based on considerations of stability, Yezzi and Mennucci have proposed a
conformal factor of the form $e^{A\ell}$ where $\ell=\int_{S^{1}}|c_{\theta
}|d\theta$, or more generally, a function $\varphi(\ell)$ of $\ell$. In this
paper, we consider only these conformal factors:
\begin{equation}
G_{c}^{\varphi}(h,k)=\varphi(\ell)\int_{S^{1}}\left(  h\cdot k\right)
|c_{\theta}|d\theta
\end{equation}
If $c(t,\theta)$ is a smooth path in $Imm(S^{1},\mathbb{R}^{2})$ connecting
$c_{o},c_{1}$, let
\begin{equation}
L_{G^{\varphi}}(c)=\int_{0}^{1}\sqrt{G_{c}^{\varphi}\left(  c_{t}^{\bot}\cdot
c_{t}^{\bot}\right)  }dt
\end{equation}
If $C_{o},C_{1}\in B_{i}(S^{1},\mathbb{R}^{2})$, define $dist_{G^{\varphi}%
}\left(  C_{1},C_{2}\right)  $ as follows. Consider their lifts $c_{o},c_{1}$
to $Imm(S^{1},\mathbb{R}^{2})$ and all paths $c(t,\theta)$ such that
$\pi(c_{o})=C_{o}$ and $\pi(c_{1})=C_{1}$. Then
\begin{equation}
dist_{G^{\varphi}}\left(  C_{1},C_{2}\right)  =\inf_{c}L_{G^{\varphi}}(c)
\end{equation}

\section{Bounds on $dist_{G^{\varphi}}$}

If $c$ is a path connecting curves $C,E$, let $\alpha(c)$ denote the area
swept out by $c$ in $\mathbb{R}^{2}$. For a path $c(t,\cdot)$, let $\ell
_{\max}(c)=\max_{t}\ell(c(t,\cdot))$. The following theorem characterizes the
$L^{1}$-type behavior of the metrics $G^{\varphi}$.

\begin{theorem}
If $\varphi(\ell)=\ell$, then,
\begin{equation}
dist_{G^{\varphi}}\left(  C,E\right)  =\inf_{c}\alpha(c)
\end{equation}

\end{theorem}

If $\varphi(\ell)=e^{A\ell}$, then,
\begin{equation}
\inf_{c}\sqrt{Ae}\alpha(c)\leq dist_{G^{\varphi}}\left(  C,E\right)  \leq
\inf_{c}\sqrt{Ae}e^{A\ell_{\max}(c)/2}\alpha(c)
\end{equation}

We first prove a series of lemmas.

\begin{lemma}%
\[
dist_{G^{\varphi}}\left(  C,E\right)  \geq\left\{
\begin{array}
[c]{l}%
\inf_{c}\alpha(c)\text{ if }\varphi(\ell)=\ell.\\
\inf_{c}\sqrt{Ae}\alpha(c)\text{ if }\varphi(\ell)=e^{A\ell}%
\end{array}
\right.
\]

\end{lemma}

\textbf{Proof: }For any path $c$,%

\begin{align*}
L_{G^{\varphi}}(c)  & =\int_{0}^{1}\left[  \varphi(\ell)\int_{S^{1}}\left(
c_{t}^{\bot}\cdot c_{t}^{\bot}\right)  |c_{\theta}|d\theta\right]
^{^{\frac{1}{2}}}dt\\
& \geq\int_{0}^{1}\left[  \left(  \frac{\varphi(\ell)}{|supp(c_{t}^{\bot}%
)|}\right)  ^{\frac{1}{2}}\int_{S^{1}}|c_{t}^{\bot}||c_{\theta}|d\theta
\right]  dt\\
& \geq\left[  \min_{t}\frac{\varphi(\ell)}{\ell}\right]  ^{\frac{1}{2}}%
\int_{S^{1}\times[0,1]}|\det dc(t,\theta)|d\theta dt\\
& \geq\left\{
\begin{array}
[c]{l}%
\alpha(c)\text{ if }\varphi(\ell)=\ell.\\
\sqrt{Ae}\alpha(c)\text{ if }\varphi(\ell)=e^{A\ell}%
\end{array}
\right.
\end{align*}
Q.E.D.

David Mumford observed from the formula for the sectional curvature that the
geodesics along which $|supp(c_{t}^{\bot})|<\ell$ if $\varphi(\ell)=\ell$ and
$<1/A$ if $\varphi(\ell)=e^{A\ell}$ may not be minimal. Such a possibility can
be heuristically seen from the inequality
\[
L_{G^{\varphi}}(c)\geq\int_{0}^{1}\left[  \left(  \frac{\varphi(\ell)}{
|supp(c_{t}^{\bot})|}\right)  ^{\frac{1}{2}}\int_{S^{1}}|c_{t}^{\bot
}||c_{\theta}|d\theta\right]  dt
\]
which suggests that while traversing a given area, one should try to minimize
$\varphi(\ell)/|supp(c_{t}^{\bot})|$. The key point is that we can increase
$|supp(c_{t}^{\bot})|$ indefinitely by replacing the part of the curve
supporting $c_{t}^{\bot}$ by a saw-tooth shaped curve of high frequency and
small amplitude. When $\varphi(\ell)=\ell$ and $|supp($ $c_{t}^{\bot})|<\ell$,
we can increase $|supp(c_{t}^{\bot})|$ so that $\ell/|supp(c_{t}^{\bot})|$
tends to $1.$ When $\varphi(\ell)=e^{A\ell}$ and $|supp(c_{t}^{\bot})|<1/A$,
we can force $e^{A|supp(c_{t}^{\bot})|}/|supp(c_{t}^{\bot})|$ to equal its
unique minimum $Ae$ by making $|supp(c_{t}^{\bot})|$ equal $1/A$. (In the case
of the metric $G^{o}$, $\varphi(\ell)=1$ so that $\varphi(\ell)/|supp(c_{t}%
^{\bot})|$ tends to $0$.) In order to obtain an upper bound for a general
path, we break it up into a series of tiny bumps. When $\varphi(\ell)=\ell$,
this method gives an upper bound for $dist_{G^{\varphi}}\left(  C,E\right)  $
which coincides with the lower bound. When $\varphi(\ell)=e^{A\ell}$, the
larger the value of $|supp(c_{t}^{\bot})|$, the greater the divergence between
the upper bound obtained by this method and the lower bound since it is more
efficient to create a large bump all at once instead a series of tiny bumps.

\textbf{Rectangular Bumps}

Let $c_{o}:S^{1}\rightarrow\mathbb{R}^{2}$ be a smooth and free immersion. Let
$C_{o}$ be the corresponding curve in $\mathbb{R}^{2}$. Let $c_{o}$ be
parametrized by the arclength so that $\theta$ parametrizes the scaled circle
$S_{\ell_{o}}^{1}$ where $\ell_{o}$ is the length of $C_{o}$. For any function
$u(\theta)$, let $u^{\prime}$ denote $du/d\theta$. Let $n_{o} $ denote the
normal vector $ic_{o}^{\prime}.$ Let $\kappa_{o}$ denote the curvature of
$c_{o}$.

Fix small positive numbers $\delta$ and $\epsilon$ such that $\delta<\ell_{o}$
and $\epsilon\left\|  \kappa_{o}\right\|  _{\infty,[0,\delta]}<<1$. Construct
a ''rectangular'' bump, over $C_{o}$ as follows:
\[
c_{1}(\theta)=\left\{
\begin{array}
[c]{ll}%
c_{o}(\theta)+\epsilon n_{o} & \text{if }0<\theta<\delta\\
\left\{  c_{o}(\theta)+sn_{o}|0\leq s\leq\epsilon\right\}  & \text{if }
\theta=0,\delta\\
c_{o}(\theta) & \text{otherwise}%
\end{array}
\right.
\]

Let $C_{1}$ be the corresponding curve in $\mathbb{R}^{2}$. The following
lemma is inspired by a comment of David Mumford that the Michor-Mumford
''teeth'' construction [2] may be used to show that the obvious path for
creating such a bump is not minimal.

\begin{lemma}
(i) If $\varphi(\ell)=\ell$,
\[
dist_{G^{\varphi}}(C_{o},C_{1})\leq\left[  \frac{1+\epsilon||\kappa
_{o}||_{\infty,[o,\delta]}}{1-\epsilon||\kappa_{o}||_{\infty,[o,\delta]}%
}\right]  ^{2}(area\;of\;the\;bump)
\]

\end{lemma}

(ii) If $\varphi(\ell)=e^{A\ell}$ and $\delta<1/A$,%

\begin{align*}
& dist_{G^{\varphi}}(C_{o},C_{1})\\
& \leq\left[  \frac{1+\epsilon||\kappa_{o}||_{\infty,[o,\delta]}}%
{1-\epsilon||\kappa_{o}||_{\infty,[o,\delta]}}\right]  ^{3/2}e^{A(\ell
_{o}+2\epsilon-\delta)/2}\sqrt{Ae^{\frac{1+\epsilon||\kappa_{o}||_{\infty
,[o,\partial]}}{1-\epsilon||\kappa_{o}||_{\infty,[o,\partial]}}}}(bump\text{
}area)
\end{align*}

\textbf{Proof: }We prove the lemma using a modification of the ''teeth''
construction of Michor and Mumford [2]. If $\varphi(\ell)=\ell$, choose
$A<1/\delta$. Approximate $C_{1}$ by a ''trapezoidal'' bump $\tilde{C}_{1}$ as
follows. Replace $C_{0}$ in the interval $[0,\delta]$ by a saw-tooth curve of
height $\eta$ and period $\frac{1}{m}$ such that its length equals $\frac
{1}{A}$. This is done by growing teeth on $C_{0}$ in time $\eta$. Move the
saw-tooth curve at unit speed along the normals $n_{0}$ keeping its end-points
fixed, until it touches the upper edge of the bump. Finally, retract the teeth
in time $\eta$. Formally, define a path $c(t,\theta)=c_{o}(\theta
)+f(t,\theta)n_{o}$ where $f(t,\theta)$ is defined as follows.
\[
f(t,\theta)=0,\text{ }0\leq t\leq\epsilon\text{ and }\delta\leq\theta\leq
\ell_{o}
\]

For $0\leq t\leq\eta,\;0\leq k\leq m-1$,
\[
f(t,\theta)=\left\{
\begin{array}
[c]{ll}%
t\left(  \frac{2m\theta}{\delta}-2k\right)  & \frac{2k}{2m}\leq\frac{\theta
}{\delta}\leq\frac{2k+1}{2m}\\
t\left(  2k+2-\frac{2m\theta}{\delta}\right)  & \frac{2k+1}{2m}\leq\frac{
\theta}{\delta}\leq\frac{2k+2}{2m}%
\end{array}
\right.
\]

For $\eta\leq t\leq\epsilon-\eta$,
\[
f(t,\theta)=\left\{
\begin{array}
[c]{ll}%
\frac{\epsilon(t-\eta)+\eta(\epsilon-\eta-t)}{\epsilon-2\eta}\cdot
\frac{2m\theta}{\delta} & 0\leq\frac{\theta}{\delta}\leq\frac{1}{2m}\\
\frac{\epsilon-\eta}{\epsilon-2\eta}(t-\eta)+f(\eta,\theta) & \frac{1}{2m}%
\leq\frac{\theta}{\delta}\leq1-\frac{1}{2m}\\
\frac{\epsilon(t-\eta)+\eta(\epsilon-\eta-t)}{\epsilon-2\eta}\cdot2m\left(
1-\frac{\theta}{\delta}\right)  & 1-\frac{1}{2m}\leq\frac{\theta}{\delta}\leq1
\end{array}
\right.
\]

For $\epsilon-\eta\leq t\leq\epsilon$,
\[
f(t,\theta)=\left\{
\begin{array}
[c]{ll}%
\frac{2m\epsilon\theta}{\delta} & 0\leq\frac{\theta}{\delta}\leq\frac{ 1}%
{2m}\\
\frac{\epsilon[t-(\epsilon-\eta)]+(\epsilon-t)f(\epsilon-\eta,\theta)}{\eta} &
\frac{1}{2m}\leq\frac{\theta}{\delta}\leq1-\frac{1}{2m}\\
2m\epsilon\left(  1-\frac{\theta}{\delta}\right)  & 1-\frac{1}{2m}\leq
\frac{\theta}{\delta}\leq1
\end{array}
\right.
\]

\begin{align*}
c^{\prime}  & =c_{o}^{\prime}+f^{\prime}n_{o}-f\kappa_{o}c_{o}^{\prime
}=(1-f\kappa_{o})c_{o}^{\prime}+f^{\prime}n_{o}\\
|c^{\prime}|  & =\sqrt{(1-f\kappa_{o})^{2}+f^{\prime2}}\\
n  & =\frac{-f^{\prime}c_{o}^{\prime}+(1-f\kappa_{o})n_{o}}{|c^{\prime}|}%
\end{align*}

\begin{align*}
c_{t}  & =f_{t}n_{o}\\
c_{t}^{\bot}\cdot c_{t}^{\bot}  & =(c_{t}\cdot n)^{2}=\frac{(1-f\kappa
_{o})^{2}f_{t}^{2}}{|c^{\prime}|^{2}}%
\end{align*}

Let
\begin{align*}
\beta & =\sqrt{\frac{(1+\epsilon||\kappa_{o}||_{\infty,[o,\delta]}%
)^{2}+f^{\prime2}}{(1-\epsilon||\kappa_{o}||_{\infty,[o,\delta]}%
)^{2}+f^{\prime2}}}\\
1  & \leq\beta\leq\frac{1+\epsilon||\kappa_{o}||_{\infty,[o,\delta]}}{
1-\epsilon||\kappa_{o}||_{\infty,[o,\delta]}}%
\end{align*}

Choose $m$ and $\eta$ such that $\int_{o}^{\delta}|c^{\prime}(\eta
,\theta)|d\theta=\frac{1}{A}$. Note that as $m\rightarrow\infty$,
$\eta\rightarrow0$.

\paragraph{Estimates when $0\leq t\leq\eta$:}

Since $|f^{\prime}|$ is independent of $\theta$ and $|f|\leq\epsilon$,
\begin{align*}
\sqrt{(1-\epsilon||\kappa_{o}||_{\infty,[o,\delta]})^{2}+f^{\prime2}}  &
\leq|c^{\prime}(\eta,\theta)|\leq\sqrt{(1+\epsilon||\kappa_{o}||_{\infty
,[o,\delta]})^{2}+f^{\prime2}}\\
\delta\sqrt{(1-\epsilon||\kappa_{o}||_{\infty,[o,\delta]})^{2}+f^{\prime2}}  &
\leq\int_{o}^{\delta}|c^{\prime}(\eta,\theta)|d\theta=\frac{1}{A}\leq
\delta\sqrt{(1+\epsilon||\kappa_{o}||_{\infty,[o,\delta]})^{2}+f^{\prime2}}%
\end{align*}

We also have%

\begin{align*}
\sqrt{(1-\epsilon||\kappa_{o}||_{\infty,[o,\delta]})^{2}+f^{\prime2}}  &
\leq|c^{\prime}(t,\theta)|\leq\sqrt{(1+\epsilon||\kappa_{o}||_{\infty
,[o,\delta]})^{2}+f^{\prime2}}\\
\frac{1}{\beta}\sqrt{(1+\epsilon||\kappa_{o}||_{\infty,[o,\delta]}%
)^{2}+f^{\prime2}}  & \leq|c^{\prime}(t,\theta)|\leq\beta\sqrt{ (1-\epsilon
||\kappa_{o}||_{\infty,[o,\delta]})^{2}+f^{\prime2}}\\
\frac{1}{\beta}\cdot\frac{1}{A\delta}  & \leq|c^{\prime}(t,\theta)|\leq
\frac{\beta}{A\delta}%
\end{align*}

\begin{align*}
\ell(c)=\int_{o}^{\ell_{o}}|c^{\prime}|d\theta & \leq\left(  \ell_{o}%
-\delta\right)  +\frac{\beta}{A}\\
e^{A\ell(c)}  & \leq e^{A\left(  \ell_{o}-\delta\right)  }e^{\beta}%
\end{align*}

Since $|f_{t}|\leq1$,
\[
\int_{o}^{\ell_{o}}|c_{t}^{\bot}|^{2}|c^{\prime}|d\theta=\int_{o}^{\delta
}\frac{(1-f\kappa_{o})^{2}f_{t}^{2}}{|c^{\prime}|} d\theta\leq\left(
1+\epsilon||\kappa_{o}||_{\infty,[o,\delta]}\right)  ^{2}\beta A\delta^{2}
\]

\[
\lim_{m\rightarrow\infty}\int_{o}^{\eta}\left[  \varphi\left(  \ell(c)\right)
\int_{o}^{\ell_{o}}|c_{t}^{\bot}|^{2}|c^{\prime}|d\theta\right]  ^{1/2}dt=0
\]

\paragraph{Estimates when $\eta\leq t\leq\epsilon-\eta$:}

Estimate for $|c^{\prime}(t,\theta)|$ is the same as in the interval
$[\frac{\delta}{2m},\delta(1-\frac{1}{2m})]$ since the curve has the same
shape. In the intervals $[0,\frac{\delta}{2m}]$ and $[\delta(1-\frac{1}{2m}
),\delta]$, $\frac{2m\eta}{\delta}\leq|f^{\prime}(t,\theta)|\leq\frac{
2m\epsilon}{\delta}$. Therefore, $|c^{\prime}(\eta,\theta)|\leq|c^{\prime
}(t,\theta)|\leq1+\epsilon||\kappa_{o}||_{\infty,[o,\partial]}+\frac
{2m\epsilon}{\delta}$.
\begin{align*}
& \ell(c)=\int_{o}^{\ell_{o}}|c^{\prime}|d\theta\leq\left(  \ell_{o}%
-\delta\right)  +\left(  1+\epsilon||\kappa_{o}||_{\infty,[o,\delta]}%
+\frac{2m\epsilon}{\delta}\right)  \frac{\delta}{m}+\frac{\beta}{A}\\
& \lim_{m\rightarrow\infty}\ell(c)\leq\left(  \ell_{o}-\delta\right)
+2\epsilon+\frac{1+\epsilon||\kappa_{o}||_{\infty,[o,\delta]}}{ 1-\epsilon
||\kappa_{o}||_{\infty,[o,\delta]}}\frac{1}{A}\\
& \lim_{m\rightarrow\infty}e^{A\ell(c)}\leq e^{A\ell_{o}+2\epsilon-\delta
)}e^{\frac{1+\epsilon||\kappa_{o}||_{\infty,[o,\delta]}}{ 1-\epsilon
||\kappa_{o}||_{\infty,[o,\delta]}}}%
\end{align*}

We also have $|f_{t}|\leq\frac{\epsilon-\eta}{\epsilon-2\eta}$. Therefore,
\[
\int_{o}^{\ell_{o}}|c_{t}^{\bot}|^{2}|c^{\prime}|d\theta\leq\left(
1+\epsilon||\kappa_{o}||_{\infty,[o,\delta]}\right)  ^{2}\left(  \frac{
\epsilon-\eta}{\epsilon-2\eta}\right)  ^{2}\beta A\delta^{2}
\]

\begin{align*}
& \lim_{m\rightarrow\infty}\int_{\eta}^{\epsilon-\eta}\left[  \varphi\left(
\ell(c)\right)  \int_{o}^{\ell_{o}}|c_{t}^{\bot}|^{2}|c^{\prime}%
|d\theta\right]  ^{1/2}dt\\
& \leq\left\{
\begin{array}
[c]{l}%
\left(  1+\epsilon||\kappa_{o}||_{\infty,[o,\delta]}\right)  \sqrt{\frac{
1+\epsilon||\kappa_{o}||_{\infty,[o,\delta]}}{1-\epsilon||\kappa_{o}%
||_{\infty,[o,\delta]}}}\left[  A\left(  \ell_{o}+2\epsilon-\delta\right)
+\frac{1+\epsilon||\kappa_{o}||_{\infty,[o,\delta]}}{1-\epsilon||\kappa
_{o}||_{\infty,[o,\delta]}}\right]  ^{1/2}\delta\epsilon\text{ }\\
\text{if }\varphi(\ell)=\ell\\
\left(  1+\epsilon||\kappa_{o}||_{\infty,[o,\delta]}\right)  \sqrt{\frac{
1+\epsilon||\kappa_{o}||_{\infty,[o,\delta]}}{1-\epsilon||\kappa_{o}%
||_{\infty,[o,\delta]}}}e^{A(\ell_{o}+2\epsilon)/2}\sqrt{Ae^{\frac{
1+\epsilon||\kappa_{o}||_{\infty,[o,\delta]}}{1-\epsilon||\kappa_{o}%
||_{\infty,[o,\delta]}}-A\delta}}\delta\epsilon\text{ }\\
\text{if }\varphi(\ell)=e^{A\ell}%
\end{array}
\right.
\end{align*}

\paragraph{Estimates when $\epsilon-\eta\leq t\leq\epsilon$:}

The path $c(t,\theta)$ in the interval $\left[  \epsilon-\eta,\epsilon\right]
$ is essentially the same as that in $\left[  0,\eta\right]  $ and
\[
\lim_{m\rightarrow\infty}\int_{\epsilon-\eta}^{\epsilon}\left[  \varphi\left(
\ell(c)\right)  \int_{o}^{\ell_{o}}|c_{t}^{\bot}|^{2}|c^{\prime}%
|d\theta\right]  ^{1/2}dt=0
\]

Since%

\[
(area\;of\;the\;bump)=\int_{o}^{\epsilon}\int_{o}^{\delta}(1-t\kappa
_{o})d\theta dt\geq\left(  1-\epsilon||\kappa_{o}||_{\infty,[o,\delta
]}\right)  \delta\epsilon\text{ }
\]
the lemma is proved in the case when $\varphi(\ell)=e^{A\ell}$. It is proved
in the case when $\varphi(\ell)=\ell$ by letting $A\rightarrow0$. Q.E.D.

In Lemma 3, we may replace the single bump by a finite number of disjoint
bumps of height $\epsilon$ and total length $\delta$. The proof remains
unchanged except that we must replace $2\epsilon$ in the formula by
$2k\epsilon$ if $k$ is the number of bumps. The function $f$ in each
individual bump may be positive or negative.

We prove the theorem by approximating the path by a series of small
rectangular bumps. The error of approximation may be made arbitrarily small by
the following lemma. Define a $Fr\acute{e}chet$ metric on $Imm(S^{1}%
,\mathbb{\ \ \ R}^{2})$ and $B_{i}(S^{1},\mathbb{R}^{2})$ as follows. If
$c_{o},c_{1}$ are points in $Imm(S^{1},\mathbb{R}^{2})$, let
\[
d_{\infty}\left(  c_{o},c_{1}\right)  =\sup_{\theta}|c_{o}\left(
\theta\right)  -c_{1}\left(  \theta\right)  |
\]
If $C_{o},C_{1}\in B_{i}(S^{1},\mathbb{R}^{2})$, let
\[
d_{\infty}\left(  C_{1},C_{2}\right)  =\inf_{\left\{  c_{o},c_{1}|\pi\left(
c_{o}\right)  =C_{o},\pi\left(  c_{1}\right)  =C_{1}\right\}  }d_{\infty
}\left(  c_{o},c_{1}\right)
\]

\begin{lemma}
For any pair
\begin{equation}
C_{1},C_{2}\in B_{i}(S^{1},\mathbb{R}^{2}),dist_{G^{\varphi}}\left(
C_{1},C_{2}\right)  \leq d_{\infty}(C_{1},C_{2})\cdot\max\{\varphi(\ell
_{1}),\varphi(\ell_{2})\}
\end{equation}
where $\ell_{i}=\ell(C_{i}),i=1,2$.
\end{lemma}

\textbf{Proof:} Let $c_{1,}c_{2}$ be lifts of $C_{1},C_{2}$ to $Imm(S^{1}%
,\mathbb{R}^{2})$. Let $c(t,\theta)=\left(  1-t\right)  c_{1}\left(
\theta\right)  +tc_{2}\left(  \theta\right)  $ be a path connecting them.
Then, $|c_{\theta}(t)|\leq\left(  1-t\right)  |c_{1,\theta}|+t|c_{2,\theta}|$
and hence, $\ell(c(t))\leq\max\{\ell(C_{1}),\ell(C_{2})\}$. Moreover,
$c_{t}=c_{2}-c_{1}.$ Therefore,
\begin{align*}
dist_{G^{\varphi}}\left(  C_{1},C_{2}\right)   & \leq\inf_{c}L_{G^{\varphi}%
}(c)\\
& \leq\inf_{\left\{  pairs\text{ }c_{1},c_{2}\right\}  }\left\{  \sup_{\theta
}\left|  c_{1}(\theta)-c_{2}(\theta)\right|  \right\}  \max\{\varphi(\ell
_{1}),\varphi(\ell_{2})\}
\end{align*}
Q.E.D.

The polygonal approximations used in the proof of the theorem lie on the
boundary of $Imm(S^{1},\mathbb{R}^{2})$ and $B_{i}(S^{1},\mathbb{R}^{2})$, and
Lemma 4 extends to them.

\textbf{Proof of the theorem:}

Consider a path $c\left(  t,\theta\right)  $ connecting $C$ and $E$. Since the
absolute curvature of the curves $c(t,\cdot)$ is uniformly bounded by a
constant $K$, each curve has a tubular neighborhood of width which is bounded
from below. Choose $\epsilon\ $ and a sequence
\[
0=t_{0}<t_{1}<\cdot\cdot\cdot<t_{N-1}<t_{N}=1
\]
such that $\epsilon K<<1$ and, for $0\leq k<N$, $c\left(  t_{k+1}%
,\cdot\right)  $ is in a local chart of $c\left(  t_{k},\cdot\right)  $:
\[
c\left(  t_{k+1},\theta\right)  =c\left(  t_{k},\theta\right)  +f_{k}\left(
\theta\right)  n_{k}
\]
where $c\left(  t_{k},\cdot\right)  $ is parametrized by the arclength, $n_{k}
$ is the normal vector field of $c\left(  t_{k},\cdot\right)  $ and
$|f_{k}|<\epsilon$. Let $F=\max\left\{  \left\|  f_{k}^{\prime}\left(
\theta\right)  \right\|  _{\infty}|0\leq k<N\right\}  $. Let $C_{k}$ denote
$\pi\left(  c\left(  t_{k},\cdot\right)  \right)  $. Let $\ell_{k}=\ell
(C_{k})$. Let $\tilde{\ell}_{k}=(1+\epsilon K+F)\ell_{k}$.

Choose $\delta$ such that
\[
\max_{k}\sqrt{\tilde{\ell}_{k}\varphi(\tilde{\ell}_{k})}\cdot F\delta
<\frac{\epsilon}{N}
\]

We now estimate the distances $dist_{G\varphi}(C_{k},$ $C_{k+1})$. Consider
the path segment $[c_{o},c_{1}]$ from $c_{o}$ to $c_{1}$ in the local chart at
$c_{o}$. Divide the range of $\theta$ into intervals of length $\delta$ .
Replace $f_{o}$ by a piecewise constant function $\bar{f}_{o}$ whose value in
each subinterval equals the average of $f_{o}$ over that interval. Let
$\bar{C}_{o}$ be the curve defined by $\bar{f}_{o}$. The $Frech\acute{e}t$
distance between $C_{1}$ and $\bar{C}_{o}$ is $\leq F\delta$. The sum of the
jumps in $\bar{f}_{o}$ is $\leq F\ell_{o}$. Since $|\bar{c}_{o}^{\prime
}|=|1-\epsilon\kappa_{o}|\leq1+\epsilon K$, $\ell(\bar{C}_{o})\leq(1+\epsilon
K+F)\ell_{o}=\tilde{\ell}_{o}$. Therefore,
\[
dist_{G^{\varphi}}\left(  C_{1},\bar{C}_{o}\right)  \leq\max\{\tilde{\ell}
_{o}\varphi(\tilde{\ell}_{o}),\tilde{\ell}_{1}\varphi(\tilde{\ell}
_{1})\}\cdot F\delta\leq\frac{\epsilon}{N}
\]

Let $\alpha([c_{o},c_{1}])$ denote the area swept out by the path $c$ during
$[0,t_{1}]$. The area between $C_{o}$ and $\bar{C}_{o}$ equals the area
between $C_{o}$ and $C_{1}$ which in turn is less than or equal to
$\alpha([c_{o},c_{1}])$. The curve $\bar{C}_{o}$ consists of a series of bumps
over $C_{o}$. Traverse the bumps sequentially, taking care to retract the
common edge of each bump with the previous bump before going to the next bump.
Retracting a common edge can be done without incurring any cost.

By Lemma 3, if $\varphi(\ell)=\ell$,
\[
dist_{G^{\varphi}}\left(  C_{o},\bar{C}_{o}\right)  \leq\left(  \frac{
1+\epsilon K}{1-\epsilon K}\right)  ^{2}\alpha([c_{o},c_{1}])
\]

and if $\varphi(\ell)=e^{A\ell}$,
\[
dist_{G^{\varphi}}\left(  C_{o},\bar{C}_{o}\right)  \leq\left(  \frac{
1+\epsilon K}{1-\epsilon K}\right)  ^{3/2}e^{A(\ell_{o}+2\epsilon)/2}\sqrt{
Ae^{\frac{1+\epsilon K}{1-\epsilon K}}}\alpha([c_{o},c_{1}])
\]

\[
dist_{G^{\varphi}}\left(  C_{o},C_{1}\right)  \leq dist_{G^{\varphi}}\left(
C_{o},\bar{C}\right)  +\frac{\epsilon}{N}
\]

Similar estimates hold for $dist_{G^{\varphi}}\left(  C_{k},C_{k+1}\right)  $
for $0<k<N$.

Therefore, if $\varphi(\ell)=\ell$,
\[
dist_{G^{\varphi}}\left(  C,E\right)  \leq\left(  \frac{1+\epsilon K}{
1-\epsilon K}\right)  ^{2}\alpha(c)+\epsilon
\]

and if $\varphi(\ell)=e^{A\ell}$,
\[
dist_{G^{\varphi}}\left(  C,E\right)  \leq\left(  \frac{1+\epsilon K}{
1-\epsilon K}\right)  ^{3/2}e^{A\ell_{\max}(c)/2}\sqrt{Ae^{\frac{1+\epsilon
K}{1-\epsilon K}}}\alpha(c)+\epsilon
\]

Since $\epsilon$ is arbitrary, we have
\[
dist_{G^{\varphi}}\left(  C,E\right)  \leq\left\{
\begin{array}
[c]{l}%
\alpha(c)\text{ if }\varphi(\ell)=\ell\\
e^{A\ell_{\max}(c)/2}\sqrt{Ae}\alpha(c)\text{ if }\varphi(\ell)=e^{A\ell}%
\end{array}
\right.
\]

Q.E.D.

For any \textit{oriented curve} $C^{or},$ define the integer-valued measurable
function $w_{C}$ on $\mathbb{R}^{2}$ by:
\[
w_{C}(x,y)=\text{the winding number of }C\text{ around }(x,y)
\]
and let
\[
d^{\flat}(C_{1}^{or},C_{2}^{or})=\int_{\mathbb{R} ^{2}}|w_{C_{1}}-w_{C_{2}%
}|dxdy
\]

It is shown in [2] that for any two \textit{oriented curves }$C_{1}^{or}%
,C_{2}^{or},$%
\[
d^{\flat}(C_{1}^{or},C_{2}^{or})\leq\min_{\text{all paths }c\text{ joining
}C_{1},C_{2}}\alpha(c)\text{ }
\]
Therefore, we have
\begin{equation}
dist_{G^{\varphi}}\left(  C_{1}^{or},C_{2}^{or}\right)  \geq\left\{
\begin{array}
[c]{l}%
d^{\flat}(C_{1}^{or},C_{2}^{or})\text{ if }\varphi(\ell)=\ell\\
\sqrt{Ae}d^{\flat}(C_{1}^{or},C_{2}^{or})\text{ if }\varphi(\ell)=e^{A\ell}%
\end{array}
\right.
\end{equation}

\begin{corollary}
(Existence of minimal geodesics) If $\varphi(\ell)=\ell$, then the only
minimal geodesics are the paths along which $|c_{t}^{\bot}|\ell$ is constant.
\end{corollary}

\textbf{Proof: }Since the inequality
\[
L_{G^{\varphi}}(c)\geq\int_{0}^{1}\left[  \int_{S^{1}}|c_{t}^{\bot}%
||c_{\theta}|d\theta\right]  dt\text{ }=\text{ }\alpha(c)
\]
is an equality if and only if $|c_{t}^{\bot}|$ does not depend on $\theta$ ,
that is, $\frac{\partial}{\partial\theta}|c_{t}^{\bot}|=0$ (the case of
''grassfire''), if $c(t)$ is a minimal geodesic, $|c_{t}^{\bot}|$ must be
independent of $\theta$. Let $c(t,\theta)$ be a path connecting $C_{1},C_{2}$
such that $|c_{t}^{\bot}|$ is independent of $\theta$. After reparametrization
if necessary, we may assume that $c_{t}\cdot c_{\theta}=0$ . Following [2], we
let $w(c)$ be the 2-current defined by the path $c(t,\theta)$. Since
$c(t,\theta)$ is an immersion,
\[
d^{\flat}(C_{1}^{or},C_{2}^{or})=\int_{\mathbb{R}^{2}}|w(c)|dxdy=\int
_{S^{1}\times[0,1]}|\det dc(t,\theta)|d\theta dt=\alpha(c)
\]
Therefore, $L_{G^{\varphi}}(c)$ is the minimal distance between $C_{1},C_{2}
$. For $c(t,\theta)$ to be a geodesic path, reparametrize $t$ such that the
infinitesimal arc-length $|c_{t}^{\bot}|\ell$ is constant along the path. Q.E.D.

\begin{corollary}
Suppose $\varphi(\ell)=e^{A\ell}$ and $c(t,\theta),0\leq t\leq1,$ is a path
connecting $C_{1},C_{2}$. Assume that $|c_{t}^{\bot}|$ does not depend on
$\theta$ and $\ell(t)<1/A$ for all $t$. Then,
\[
dist_{G^{\varphi}}(C_{1}^{or},C_{2}^{or})=\sqrt{Ae}d^{\flat}(C_{1}^{or}%
,C_{2}^{or})
\]
It follows that the path $c(t,\theta)$ is not minimal.
\end{corollary}

\textbf{Proof: }Break up the interval $[0,1]$ into small segments of length
$\epsilon$ and apply Lemma 3 with $\delta=\ell_{o}$. (The proof of the lemma
extends to this case after minor modifications.) Calculate the length of the
new path $\tilde{c}_{\epsilon}$ applying the construction of Lemma 3 to each
of the segments.$.$ We get $\lim_{_{\epsilon\rightarrow0}}L_{G^{\varphi}%
}(\tilde{c}_{\epsilon})=\sqrt{Ae}\alpha(c)$ as in the proof of Theorem 2. On
the other hand,
\begin{align*}
L_{G^{\varphi}}(c)  & =\int_{0}^{1}\left[  \int_{S^{1}}|c_{t}^{\bot
}||c_{\theta}|d\theta\right]  dt=\int_{o}^{1}\sqrt{\frac{\varphi(\ell)}{ \ell
}}|c_{t}^{\bot}|\ell dt\\
& >\sqrt{\min_{\ell}\frac{\varphi(\ell)}{\ell}}\alpha(c)=\sqrt{Ae}\alpha(c)
\end{align*}
Q.E.D.

\section{Geodesic Equations}

We reproduce calculations in [2] mutatis mutandis. Let $t\longmapsto
c(t,\cdot)$ be a path in $Imm(S^{1},\mathbb{R}^{2})$.

\textbf{Equation of Geodesic in }$Imm(S^{1},\mathbb{R}^{2}):$%
\begin{equation}
\left(  \varphi(\ell)|c_{\theta}|c_{t}\right)  _{t}=-\frac{1}{2}\left[
\left(  \int_{S^{1}}\left|  c_{t}\right|  ^{2}\left|  c_{\theta}\right|
d\theta\right)  i\varphi^{\prime}(\ell)\kappa_{c}c_{\theta}+\varphi
(\ell)\left(  \frac{\left|  c_{t}\right|  ^{2}c_{\theta}}{\left|  c_{\theta
}\right|  }\right)  _{\theta}\right]
\end{equation}
where $\varphi^{\prime}$ denotes $d\varphi/d\ell$.

\textbf{Proof: }We calculate the first variation of the energy of the path to
obtain the geodesics:.%

\[
E_{G^{\varphi}}\left(  c\right)  =\frac{1}{2}\int_{a}^{b}\left[  \varphi
(\ell)\left(  \int_{S^{1}}\left(  c_{t}\cdot c_{t}\right)  |c_{\theta}%
|d\theta\right)  \right]  dt
\]
\begin{align*}
\partial_{s}|_{o}\varphi(\ell)  & =\varphi^{\prime}(\ell)\partial_{s}|_{o}%
\int_{S^{1}}\left(  c_{\theta}\cdot c_{\theta}\right)  ^{\frac{1}{2} }%
d\theta\\
& =\varphi^{\prime}(\ell)\int_{S^{1}}\frac{\left(  c_{s\theta}\cdot c_{\theta
}\right)  }{\left|  c_{\theta}\right|  }d\theta\\
& =-\varphi^{\prime}(\ell)\int_{S^{1}}(c_{s}\cdot\left(  \frac{c_{\theta}%
}{\left|  c_{\theta}\right|  }\right)  _{\theta})d\theta\\
& =-\varphi^{\prime}(\ell)\int_{S^{1}}\left(  c_{s}\cdot i\kappa_{c}c_{\theta
}\right)  d\theta
\end{align*}
\begin{align*}
\partial_{s}|_{o}\frac{1}{2}\int_{S^{1}}\left(  c_{t}\cdot c_{t}\right)
|c_{\theta}|d\theta & =\int_{S^{1}}\left(  c_{st}\cdot c_{t}\right)
|c_{\theta}|d\theta+\frac{1}{2}\int_{S^{1}}\left|  c_{t}\right|  ^{2}%
\frac{\left(  c_{s\theta}\cdot c_{\theta}\right)  }{\left|  c_{\theta}\right|
} d\theta\\
& =\int_{S^{1}}\left(  c_{st}\cdot c_{t}\right)  |c_{\theta}|d\theta-\frac
{1}{2}\int_{S^{1}}(c_{s}\cdot\left(  \frac{\left|  c_{t}\right|  ^{2}%
c_{\theta}}{\left|  c_{\theta}\right|  }\right)  _{\theta})d\theta
\end{align*}
\begin{align*}
\partial_{s}|_{o}E_{G^{\varphi}}\left(  c\right)   & =-\frac{1}{2}%
\varphi^{\prime}(\ell)\int_{a}^{b}\left(  \int_{S^{1}}\left(  c_{s}\cdot
i\kappa_{c}c_{\theta}\right)  d\theta\right)  \left(  \int_{S^{1}}\left|
c_{t}\right|  ^{2}|c_{\theta}|d\theta\right)  dt\\
& -\frac{1}{2}\int_{a}^{b}\varphi(\ell)\left(  \int_{S^{1}}(c_{s}\cdot\left(
\frac{\left|  c_{t}\right|  ^{2}c_{\theta}}{\left|  c_{\theta}\right|
}\right)  _{\theta})d\theta\right)  dt\\
& +\int_{a}^{b}\varphi(\ell)\left(  \int_{S^{1}}\left(  c_{st}\cdot
c_{t}\right)  |c_{\theta}|d\theta\right)  dt
\end{align*}

\[
The\;last\;term\;=-\int_{a}^{b}\int_{S^{1}}\left(  c_{s}\cdot\varphi
(\ell)c_{t}\left|  c_{\theta}\right|  _{t}\right)  d\theta dt
\]
Therefore,
\begin{align*}
\partial_{s}|_{o}E_{G^{\varphi}}\left(  c\right)   & =-\int_{a}^{b}\int
_{S^{1}}\left(  c_{s}\cdot F\right)  d\theta dt\\
\text{where}\quad F  & =\left(  \varphi(\ell)\left|  c_{\theta}\right|
c_{t}\right)  _{t}+\frac{i\varphi^{\prime}(\ell)\kappa_{c}c_{\theta}}%
{2}\left(  \int_{S^{1}}\left|  c_{t}\right|  ^{2}|c_{\theta}|d\theta\right) \\
& +\frac{\varphi(\ell)}{2}\left(  \frac{\left|  c_{t}\right|  ^{2}c_{\theta}%
}{\left|  c_{\theta}\right|  }\right)  _{\theta}%
\end{align*}
Q.E.D.

\textbf{Equation of Geodesic in }$B_{i}(S^{1},\mathbb{R}^{2}):$ To obtain the
geodesics in $B_{i}(S^{1},\mathbb{R}^{2}),$ we may assume that $\left(
c_{t}\cdot c_{\theta}\right)  =0$. Write $c_{t}$ as $iac_{\theta}/\left|
c_{\theta}\right|  $. After substituting this in Eq. (9), split the equation
into its component along $c_{\theta}$ and $n_{c}.$ The former vanishes
identically. The normal component of the left-hand-side is
\[
\left(  \varphi(\ell)|c_{\theta}|a\right)  _{t}=\varphi(\ell)a_{t}\left|
c_{\theta}\right|  -\varphi(\ell)\kappa_{c}a^{2}\left|  c_{\theta}\right|
-\varphi^{\prime}(\ell)a\left|  c_{\theta}\right|  \int_{S^{1}}a\kappa
_{c}\left|  c_{\theta}\right|  d\theta
\]
where we have used the formulae $\left|  c_{\theta}\right|  _{t}=-a\kappa
_{c}\left|  c_{\theta}\right|  $ and $\ell_{t}(c)=-\int_{S^{1}}a\kappa
_{c}\left|  c_{\theta}\right|  d\theta$. The normal component of the right
hand side is
\[
-\frac{1}{2}\varphi(\ell)\kappa_{c}a^{2}\left|  c_{\theta}\right|  -\frac
{1}{2}\varphi^{\prime}(\ell)\kappa_{c}\left(  \int_{S^{1}}a^{2}\left|
c_{\theta}\right|  d\theta\right)  \left|  c_{\theta}\right|
\]

Therefore, the equation of the geodesic in $B_{i}(S^{1},\mathbb{R}^{2})$ is
\begin{align}
a_{t}  & =\frac{\kappa_{c}}{2}\left(  a^{2}-\frac{\varphi^{\prime}(\ell)}{
\varphi(\ell)}\int_{S^{1}}a^{2}\left|  c_{\theta}\right|  d\theta\right)
+\frac{\varphi^{\prime}(\ell)}{\varphi(\ell)}a\int_{S^{1}}a\kappa_{c}\left|
c_{\theta}\right|  d\theta\nonumber\\
& =\left\{
\begin{array}
[c]{ll}%
\frac{\kappa_{c}}{2}\left(  a^{2}-\overline{a^{2}}\right)  +a\cdot\overline{
a\kappa_{c}} & \text{ if }\varphi(\ell)=\ell\\
\frac{\kappa_{c}}{2}\left(  a^{2}-(A\ell)\overline{a^{2}}\right)
+(A\ell)a\cdot\overline{a\kappa_{c}} & \text{if }\varphi(\ell)=e^{A\ell}%
\end{array}
\right.
\end{align}
where for any function $f(t,\theta)$, $\overline{f}$ denotes the average
$\frac{1}{\ell}\int_{S^{1}}f\left|  c_{\theta}\right|  d\theta$.

As an example, consider the case of concentric circles [2], $c(t,\theta
)=r(t)e^{i\theta},r_{o}=r(0),r_{1}=r(1)$ We have $\kappa_{c}=1/r$ and
$a=-r_{t}.$ Substituting these in Eq. (10) when $\varphi(\ell)=\ell$, we get
$-r_{tt}=r_{t}^{2}/r$ or $\left(  r^{2}\right)  _{tt}=0$. Therefore,
$r^{2}\left(  t\right)  =tr_{1}^{2}+(1-t)r_{o}^{2}$. This example is a special
case of the curve evolution by ''grassfire'' in which $a$ is independent of
$\theta$. We have $a^{2}=\overline{a^{2}}$ and the equation of the geodesic
reduces to $a_{t}=a\cdot\overline{a\kappa_{c}}=-a\ell_{t}/\ell$ and hence
$\left(  a\ell\right)  _{t}=0$. Therefore, $a\ell=$ a constant. By
substituting in the equation for the length of the geodesic, we find that
$a\ell$ = the length of the geodesic. When $\varphi(\ell)=e^{A\ell}$, the
equation of the geodesic in the case of concentric circles is $\left(
r^{2}\right)  _{tt}=r_{t}^{2}(1-2\pi rA)$ which is zero when the perimeter of
the circle equals $1/A$, marking the unique inflection point of the function
$r^{2}(t)$.

\section{Sectional Curvature}

A formula for the sectional curvature may be derived exactly as in [2] by
means of local charts. Let $c:S^{1}\rightarrow\mathbb{R}^{2}$ be a smooth and
free immersion, an element in the space of free immersions, $Imm_{f}%
(S^{1},\mathbb{R}^{2})$. $B_{i,f}(S^{1},\mathbb{R}^{2})=\pi\left(
Imm_{f}(S^{1},\mathbb{R} ^{2})\right)  $. Let $C=$ $\pi(c)$. As before, let
$c$ be parametrized by the arclength so that $\theta$ parametrizes the scaled
circle $S_{\ell}^{1} $ where $\ell$ is the length of $c$. Let
\begin{align*}
\psi:C^{\infty}(S_{\ell}^{1},(-\epsilon,\epsilon))  & \rightarrow
Imm_{f}(S^{1},\mathbb{R}^{2})\\
\psi(f)(\theta)  & =c(\theta)+f(\theta)n_{c}(\theta)\\
\pi\circ\psi:C^{\infty}(S_{\ell}^{1},(-\epsilon,\epsilon))  & \rightarrow
B_{i,f}(S^{1},\mathbb{R}^{2})
\end{align*}
be a local chart centered at $C$. For any function $u$ , let $u^{\prime}$
denote its derivative so that $\varphi^{\prime}=d\varphi/d\ell$ and if $u $ is
a function on $S_{\ell}^{1}$, $u^{\prime}=du/d\theta$; . We have the following
formulae from [2]:
\begin{align*}
\psi(f)^{\prime}  & =(1-f\kappa_{c})c^{\prime}+f^{\prime}n_{c}\\
n_{\psi(f)}  & =\frac{(1-f\kappa_{c})n_{c}-f^{\prime}c^{\prime}}%
{\sqrt{(1-f\kappa_{c})^{2}+f^{\prime2}}}%
\end{align*}
For $h\in C^{\infty}(S_{\ell}^{1},\mathbb{R}),$ $h\cdot n_{c}\in T_{\psi
(f)}Imm_{f}(S^{1},\mathbb{R}^{2})$ and
\[
\left(  h\cdot n_{c}\right)  ^{\bot}=\frac{(1-f\kappa_{c})h}{\sqrt
{(1-f\kappa_{c})^{2}+f^{\prime2}}}n_{\psi(f)}
\]
Let
\begin{align*}
<h,k>_{f}  & =\int_{S_{\ell}^{1}}\frac{(1-f\kappa_{c})^{2}hk}{\sqrt{
(1-f\kappa_{c})^{2}+f^{\prime2}}}d\theta\\
& =\int_{S_{\ell}^{1}}hk(1-f\kappa_{c}-\frac{1}{2}f^{\prime2}+O\left(
f^{3}\right)  )d\theta
\end{align*}
\[
\ell_{f}=\int_{S_{\ell}^{1}}\sqrt{(1-f\kappa_{c})^{2}+f^{\prime2}}
d\theta=\int_{S_{\ell}^{1}}(1-f\kappa_{c}+\frac{1}{2}f^{\prime2}+O\left(
f^{3}\right)  )d\theta
\]
$<h,k>_{o}$ will be denoted simply as $<h,k>$ and $||h||^{2}=<h,h>$. The
metric at the point in $B_{i,f}(S^{1},\mathbb{R}^{2})$ corresponding to $f$ is
given by
\[
G_{f}^{\varphi}(h,k)=\varphi(\ell_{f})<h,k>_{f}
\]
We now outline the calculation of the Christoffel Symbol and the sectional
curvature. (For more details, see [2].) Differentiating the metric in
direction $j$, we get
\[
d_{j}G_{f}^{\varphi}(h,k)=\left(  \varphi^{\prime}(\ell_{f})d_{j}\ell
_{f}\right)  <h,k>_{f}+\varphi(\ell_{f})d_{j}<h,k>_{f}
\]

\[
d_{j}\ell_{f}=-<j,\kappa_{c}>+<f^{\prime},j^{\prime}>+O(f^{2})\,=-<\kappa
_{c}+f^{\prime\prime},j>+O(f^{2})
\]

\[
d_{j}<h,k>_{f}=-<hk,j\kappa_{c}+f^{\prime}j^{\prime}>+O(f^{2})
\]
Therefore,
\begin{align*}
& d_{j}G_{f}^{\varphi}(h,k)\\
& =\varphi^{\prime}(\ell_{f})\left(  -<j,\kappa_{c}>+<f^{\prime},j^{\prime
}>\right)  \,<h,k>_{f}-\,\varphi(\ell_{f})<hk,j\kappa_{c}+f^{\prime}j^{\prime
}>+O(f^{2})\\
& =-\varphi^{\prime}(\ell_{f})<\kappa_{c}+f^{\prime\prime},j>\,<h,k>_{f}%
-\,\varphi(\ell_{f})<hk,j\kappa_{c}+f^{\prime}j^{\prime}>+O(f^{2})
\end{align*}

The Christoffel symbol satisfies the identity
\[
2G_{f}^{\varphi}\left(  \Gamma_{f}\left(  h,k\right)  ,j\right)  =d_{j}%
G_{f}^{\varphi}(h,k)-\;d_{h}G_{f}^{\varphi}\left(  k,j\right)  -\,d_{k}%
G_{f}^{\varphi}\left(  h,j\right)
\]
Therefore,
\begin{align*}
2G_{f}^{\varphi}\left(  \Gamma_{f}\left(  h,k\right)  ,j\right)
=\varphi^{\prime}(\ell_{f})<\kappa_{c}+f^{\prime\prime},h>\,<k,j>_{f}%
+\,\varphi(\ell_{f})<kj,h\kappa_{c}+f^{\prime}h^{\prime}>  & \\
+\,\varphi^{\prime}(\ell_{f})<\kappa_{c}+f^{\prime\prime},k>\,<h,j>_{f}%
+\,\varphi(\ell_{f})<hj,k\kappa_{c}+f^{\prime}k^{\prime}>  & \\
-\varphi^{\prime}(\ell_{f})\,<\kappa_{c}+f^{\prime\prime},j>\,<h,k>_{f}%
-\,\varphi(\ell_{f})<hk,j\kappa_{c}+f^{\prime}j^{\prime}>+O(f^{2})  &
\end{align*}
At the center $f=0$, this simplifies to
\[
2G_{o}^{\varphi}\left(  \Gamma_{o}\left(  h,k\right)  ,j\right)
=\,<\varphi(\ell)\kappa_{c}hk+\varphi^{\prime}(\ell)\left[  <\kappa
_{c},h>k+<\kappa_{c},k>h\,-\,<h,k>\kappa_{c}\right]  ,\,j>
\]
Since $G_{o}^{\varphi}\left(  \Gamma_{o}\left(  h,k\right)  ,j\right)
=\varphi(\ell)<\Gamma_{o}\left(  h,k\right)  ,j>$, we have
\[
\Gamma_{o}\left(  h,k\right)  =\frac{1}{2}\kappa_{c}hk+\frac{\varphi^{\prime
}(\ell)}{2\varphi(\ell)}\left[  <\kappa_{c},h>k+<\kappa_{c}%
,k>h\,-\,<h,k>\kappa_{c}\right]
\]

Next, we calculate the second derivative:
\begin{align*}
& d_{m}d_{j}G_{f}^{\varphi}(h,k)\\
& =d_{m}\left(  \varphi^{\prime}(\ell_{f})\left(  -<j,\kappa_{c}>+<f^{\prime
},j^{\prime}>\right)  \,<h,k>_{f}-\,\varphi(\ell_{f})<hk,j\kappa_{c}%
+f^{\prime}j^{\prime}>+O(f^{2})\right) \\
& =\varphi^{\prime\prime}(\ell_{f})\left(  <m,\kappa_{c}>-<f^{\prime
},m^{\prime}>\right)  \left(  <j,\kappa_{c}>-<f^{\prime},j^{\prime}>\right)
<h,k>_{f}\\
& \;+\,\varphi^{\prime}(\ell_{f})<m^{\prime},j^{\prime}><h,k>_{f}%
+\varphi^{\prime}(\ell_{f})\left(  <j,\kappa_{c}>-<f^{\prime},j^{\prime
}>\right)  <hk,m\kappa_{c}+f^{\prime}m^{\prime}>\\
& \,+\varphi^{\prime}(\ell_{f})\left(  <m,\kappa_{c}>-<f^{\prime},m^{\prime
}>\right)  <hk,j\kappa_{c}+f^{\prime}j^{\prime}>-\,\,\varphi(\ell
_{f})<hk,m^{\prime}j^{\prime}>+O\left(  f\right)
\end{align*}

At the center,
\begin{align*}
& d_{m}d_{j}G_{o}^{\varphi}(h,k)\\
& =\varphi^{\prime\prime}(\ell)<m,\kappa_{c}>\,<j,\kappa_{c}>\,<h,k>\,\\
& +\varphi^{\prime}(\ell)\left(  <m^{\prime},j^{\prime}>\,<h,k>+<j,\kappa
_{c}>\,<hk,m\kappa_{c}>+<m,\kappa_{c}>\,<hk,j\kappa_{c}>\right) \\
& \;-\,\,\varphi(\ell)<hk,m^{\prime}j^{\prime}>
\end{align*}

Let $P(m,h)$ be a tangent plane at $c,$ spanned by normal vector fields
$m(\theta),h\left(  \theta\right)  $ along $C$. Assume that $m,h$ have been
normalized so that
\[
\varphi(\ell)\left\|  m\right\|  ^{2}=\varphi(\ell)\left\|  h\right\|
^{2}=1\text{ and }<m,h>=0
\]
The sectional curvature at $C$ is given by the formula
\begin{align*}
k_{c}(P(m,h))  & =d_{m}d_{h}G_{o}^{\varphi}(h,m)-\frac{1}{2}\left[  d_{m}%
d_{m}G_{o}^{\varphi}(h,h)+d_{h}d_{h}G_{o}^{\varphi}(m,m)\right] \\
& +G_{o}^{\varphi}\left(  \Gamma_{o}\left(  h,m\right)  ,\Gamma_{o}\left(
h,m\right)  \right)  -G_{o}^{\varphi}\left(  \Gamma_{o}\left(  m,m\right)
,\Gamma_{o}\left(  h,h\right)  \right)
\end{align*}

\begin{align*}
& d_{m}d_{h}G_{o}^{\varphi}(h,m)-\frac{1}{2}\left[  d_{m}d_{m}G_{o}^{\varphi
}(h,h)+d_{h}d_{h}G_{o}^{\varphi}(m,m)\right] \\
& =\varphi^{\prime}\,\left(  <m,\kappa_{c}>\,<mh,h\kappa_{c}>+<h,\kappa
_{c}>\,<mh,m\kappa_{c}>\right)  -\,\,\varphi(\ell)<mh,m^{\prime}h^{\prime}>\\
& -\frac{\varphi^{\prime\prime}}{2}\left[  <m,\kappa_{c}>^{2}\,\left\|
h\right\|  ^{2}+<h,\kappa_{c}>^{2}\,\left\|  m\right\|  ^{2}\right] \\
& -\frac{\varphi^{\prime}}{2}\left[  \,\left\|  m^{\prime}\right\|
^{2}\,\left\|  h\right\|  ^{2}+2<m,\kappa_{c}>\,<h^{2},m\kappa_{c}>\right]
+\frac{\varphi}{2}<m^{\prime2},h^{2}>\\
& \,-\frac{\varphi^{\prime}}{2}\left[  \left\|  m\right\|  ^{2}\,\left\|
h^{\prime}\right\|  ^{2}+2<h,\kappa_{c}>\,<m^{2},h\kappa_{c}>\right]
+\frac{\varphi}{2}\,<m^{2},h^{\prime2}>\\
& =\frac{\varphi}{2}\left\|  m^{\prime}h-mh^{\prime}\right\|  ^{2}-\frac{
\varphi^{\prime}}{2\varphi}\left(  \left\|  m^{\prime}\right\|  ^{2}+\left\|
h^{\prime}\right\|  ^{2}\right)  -\frac{\varphi^{\prime\prime}}{2\varphi
}\left(  <m,\kappa_{c}>^{2}+<h,\kappa_{c}>^{2}\right)
\end{align*}

\begin{align*}
& 4\left[  \Gamma_{o}^{2}\left(  h,m\right)  -\Gamma_{o}\left(  m,m\right)
\Gamma_{o}\left(  h,h\right)  \right] \\
& =\left[  \frac{\varphi^{\prime}}{\varphi}\left(  <m,\kappa_{c}%
>h+<h,\kappa_{c}>m\right)  +mh\kappa_{c}\right]  ^{2}\\
& \,-\left[  \frac{\varphi^{\prime}}{\varphi}\left(  2<m,\kappa_{c}%
>m-\kappa_{c}\,\left\|  m\right\|  ^{2}\right)  +m^{2}\kappa_{c}\right]
\left[  \frac{\varphi^{\prime}}{\varphi}\left(  2<h,\kappa_{c}>h-\kappa
_{c}\,\left\|  h\right\|  ^{2}\right)  +h^{2}\kappa_{c}\right] \\
& =\left(  \frac{\varphi^{\prime}}{\varphi}\right)  ^{2}\left(  <m,\kappa
_{c}>^{2}h^{2}+<h,\kappa_{c}>^{2}m^{2}\right) \\
& +\left(  \frac{\varphi^{\prime}}{\varphi}\right)  ^{2}\left(  2<m,\kappa
_{c}>m\kappa_{c}\left\|  h\right\|  ^{2}+2<h,\kappa_{c}>h\kappa_{c}\left\|
m\right\|  ^{2}-\kappa_{c}^{2}\left\|  h\right\|  ^{2}\left\|  m\right\|
^{2}\right) \\
& -2\left(  \frac{\varphi^{\prime}}{\varphi}\right)  ^{2}<m,\kappa
_{c}>\,<h,\kappa_{c}>mh+\frac{\varphi^{\prime}}{\varphi}\left(  m^{2}%
\kappa_{c}^{2}\left\|  h\right\|  ^{2}+h^{2}\kappa_{c}^{2}\left\|  m\right\|
^{2}\right)
\end{align*}

Taking into account that $h,m$ form an orthonormal basis,
\begin{align*}
& G_{o}^{\varphi}\left(  \Gamma_{o}\left(  h,m\right)  ,\Gamma_{o}\left(
h,m\right)  \right)  -G_{o}^{\varphi}\left(  \Gamma_{o}\left(  m,m\right)
,\Gamma_{o}\left(  h,h\right)  \right) \\
& =\frac{1}{4}\left(  \frac{\varphi^{\prime}}{\varphi}\right)  ^{2}\left(
3<m,\kappa_{c}>^{2}+\,3<h,\kappa_{c}>^{2}-\frac{\left\|  \kappa_{c}\right\|
^{2}}{\varphi}\right)  +\frac{1}{4}\frac{\varphi^{\prime}}{ \varphi}\left(
\left\|  m\kappa_{c}\right\|  ^{2}+\left\|  h\kappa_{c}\right\|  ^{2}\right)
\end{align*}

Putting all of this together, we get the sectional curvature%

\begin{align}
& k_{c}(P(m,h))=\frac{\varphi}{2}\left\|  m^{\prime}h-mh^{\prime}\right\|
^{2}-\frac{\varphi^{\prime}}{2\varphi}\left(  \left\|  m^{\prime}\right\|
^{2}+\left\|  h^{\prime}\right\|  ^{2}\right)  +\frac{\varphi^{\prime}}{
4\varphi}\left(  \left\|  m\kappa_{c}\right\|  ^{2}+||h\kappa_{c}||^{2}\right)
\nonumber\\
& +\frac{3\varphi^{\prime2}-2\varphi\varphi^{\prime\prime}}{4\varphi^{2}%
}\left(  <m,\kappa_{c}>^{2}+\,<h,\kappa_{c}>^{2}\right)  -\frac{\varphi
^{\prime2}}{4\varphi^{3}}\left\|  \kappa_{c}\right\|  ^{2}%
\end{align}

Each of the last three terms on the right-hand side is bounded:
\[
\left\|  m\kappa_{c}\right\|  ^{2}+\left\|  h\kappa_{c}\right\|  ^{2}\leq
\frac{2\left\|  \kappa_{c}\right\|  _{\infty}^{2}}{\varphi},\;<m,\kappa
_{c}>^{2}+\,<h,\kappa_{c}>^{2}\leq\frac{2\ell\left\|  \kappa_{c}\right\|
_{\infty}^{2}}{\varphi},\;\left\|  \kappa_{c}\right\|  ^{2}\leq\ell\left\|
\kappa_{c}\right\|  _{\infty}^{2}
\]
where $\left\|  \kappa_{c}\right\|  _{\infty}=\max_{\theta}|\kappa_{c}%
(\theta)|$. Therefore, the boundedness of the sectional curvature from above
and hence the minimality of a geodesic crucially depend on the first two
terms. For a fixed $m$, the magnitude of each of the two terms depends on
$||h^{^{\prime}}||$ which can be made arbitrarily large while keeping $||h||$
fixed by making $h$ highly wiggly.

\begin{proposition}
For a given $m$, the sectional curvature is bounded from above if and only if
$||m||_{\infty}^{2}\leq\dfrac{\varphi^{\prime}}{\varphi^{2}}$ or,
equivalently, $m^{2}(\theta)\leq$.$\dfrac{\varphi^{\prime}\ell}{\varphi
}\overline{m^{2}}$ since $\dfrac{\varphi^{\prime}}{\varphi^{2}}=\dfrac
{\varphi^{\prime}}{\varphi^{2}}\cdot\varphi\int_{S_{\ell}^{1}}m^{2}d\theta$.
\end{proposition}

\textbf{Proof: }We need to estimate only the first two terms on the right-hand
side of Eq. (11).

Suppose $||m||_{\infty}^{2}\leq\dfrac{\varphi^{\prime}}{\varphi^{2}}$. Then,
\begin{align*}
& \frac{\varphi}{2}\left\|  m^{\prime}h-mh^{\prime}\right\|  ^{2}%
-\frac{\varphi^{\prime}}{2\varphi}\left(  \left\|  m^{\prime}\right\|
^{2}+\left\|  h^{\prime}\right\|  ^{2}\right) \\
& =\varphi\int_{S_{\ell}^{1}}\left[  \frac{\varphi}{2}(m^{\prime}h)^{2}%
-\frac{\varphi^{\prime}}{2\varphi}m^{\prime2}-\varphi mm^{\prime}hh^{\prime
}+\frac{\varphi}{2}\left(  m^{2}-\frac{\varphi^{\prime}}{\varphi^{2}}\right)
h^{\prime2}\right]  d\theta\\
& \leq\frac{\varphi}{2}\left\|  m^{\prime}\right\|  _{\infty}^{2}%
-\frac{\varphi^{2}}{4}\int_{S_{\ell}^{1}}(m^{2})^{\prime}(h^{2})^{\prime
}d\theta\\
& \leq\frac{\varphi}{2}\left\|  m^{\prime}\right\|  _{\infty}^{2}%
+\frac{\varphi^{2}}{4}\int_{S_{\ell}^{1}}(m^{2})^{\prime\prime}(h^{2}%
)d\theta\\
& \leq\frac{\varphi}{2}\left\|  m^{\prime}\right\|  _{\infty}^{2}+\frac{
\varphi}{4}\left\|  (m^{2})^{\prime\prime}\right\|  _{\infty}<\infty
\end{align*}

Conversely, suppose $||m||_{\infty}^{2}>\dfrac{\varphi^{\prime}}{\varphi^{2}}%
$. Choose $\epsilon$ such that $U=\{\theta:m^{2}(\theta)>\varphi^{\prime
}/\varphi^{2}+\epsilon\}$ is not empty. Let $h$ be a high frequency wave
function with $supp(h)\subset U$. Then,
\begin{align*}
& \frac{\varphi}{2}\left\|  m^{\prime}h-mh^{\prime}\right\|  ^{2}-\frac{
\varphi^{\prime}}{2\varphi}\left(  \left\|  m^{\prime}\right\|  ^{2}+\left\|
h^{\prime}\right\|  ^{2}\right) \\
& =\varphi\int_{S_{\ell}^{1}}\left[  \frac{\varphi}{2}(m^{\prime}h)^{2}%
-\frac{\varphi^{\prime}}{2\varphi}m^{\prime2}-\varphi mm^{\prime}hh^{\prime
}+\frac{\varphi}{2}\left(  m^{2}-\frac{\varphi^{\prime}}{ \varphi^{2}}\right)
h^{\prime2}\right]  d\theta\\
& \geq-\frac{\varphi^{\prime}}{2\varphi}||m^{\prime}||^{2}+\frac{\varphi^{2}%
}{4}\int_{S_{\ell}^{1}}(m^{2})^{\prime\prime}(h^{2})d\theta+\frac
{\varphi\epsilon}{2}\left\|  h^{\prime}\right\|  ^{2}\\
& \geq-\frac{\varphi^{\prime}}{2\varphi}||m^{\prime}||^{2}-\frac{\varphi^{2}%
}{4}\left\|  (m^{2})^{\prime\prime}\right\|  _{\infty}+\frac{ \varphi\epsilon
}{2}\left\|  h^{\prime}\right\|  ^{2}%
\end{align*}
which tends to $\infty$ as the frequency of the wave function $h\ $tends to
$\infty$. Q.E.D.

If $U=\{\theta:m^{2}(\theta)>\varphi^{\prime}/\varphi^{2}+\epsilon\}$ is not
empty,
\begin{align*}
1  & =\varphi\left[  \int_{U}m^{2}d\theta+\int_{[0,\ell]/U}m^{2}d\theta\right]
\\
& \geq\varphi\int_{U}m^{2}d\theta\geq\left(  \frac{\varphi^{\prime}}{\varphi
}+\epsilon\varphi\right)  |U|
\end{align*}
and hence, $|U|<\varphi/\varphi^{\prime}$. If $\varphi=\ell$, $|U|<\ell$ and
if $\varphi=e^{A\ell}$, $|U|<1/A$. Thus, the case when $||m||_{\infty}%
^{2}>\varphi^{\prime}/\varphi^{2}$ may be seen as a generalization of the
rectangular bump considered in \S 3.

If $\varphi=\ell$, the sectional curvature is bounded if and only if
$m^{2}(\theta)\leq$.$\overline{m^{2}}$ which is true if and only if $m=1/\ell
$. Setting $m=1/\ell$, we get%

\[
k_{c}(P(\frac{1}{\ell},h))=\frac{\left\|  h\kappa_{c}\right\|  ^{2}}{4\ell
}+\frac{3}{4\ell^{2}}\left(  <\frac{1}{\ell},\kappa_{c}>^{2}+\,<h,\kappa
_{c}>^{2}\right)
\]
which is always positive. If $h$ and $\kappa_{c}$ additionally have disjoint
supports, the sectional curvature equals $3\pi^{2}N^{2}/\ell^{4}$ where $N$ is
the rotation index of $C$.

If $\varphi=e^{A\ell}$, the sectional curvature is bounded if and only if
$m^{2}(\theta)\leq$.$(A\ell)\overline{m^{2}}$. In particular, the sectional
curvature is unbounded for every $m$ if $\ell<1/A$. The analysis of
rectangular bumps in \S 3 suggests the conjecture that when $\varphi=e^{A\ell
}$, a geodesic is locally minimal if and only if $a^{2}(\theta)\leq
(A\ell)\overline{a^{2}}$ where $a=|c_{t}^{\bot}|$ (notation of Eq. 10).

For an example of a negative sectional curvature, consider the unit square
with slightly rounded corners. Choose $m,h$ such that $supp(m)$ and $supp(h)$
are disjoint and concentrated along the straight portions of the square.
Then,
\[
k_{c}(P(m,h))=-\frac{\varphi^{\prime}}{2\varphi}\left(  \left\Vert m^{\prime
}\right\Vert ^{2}+\left\Vert h^{\prime}\right\Vert ^{2}\right)  -\frac
{\varphi^{\prime2}}{4\varphi^{3}}\left\Vert \kappa_{c}\right\Vert ^{2}
\]

\textbf{Acknowledgement: }Suggestions of David Mumford contributed greatly to
this paper.

\section{References}

\begin{enumerate}
\item G. Charpiat, O. Faugeras and R. Keriven, ''Approximations of shape
metrics and application to shape warping and empirical shape statistics'', To
appear in J. of Foundations of Computational Mathematics.

\item P. Michor and D. Mumford, ''Riemannian geometries on spaces of plane
curves'', arXiv:math.DG/0312384, v2, Sep. 22, 2004.

\item E. Klassen, A. Srivastava, W. Mio and S.H. Joshi, ''Analysis of planar
shapes using geodesic paths on shape spaces'', IEEE Trans. PAMI, 26(3), pp.
372-383, 2004.

\item A. Yezzi and A. Mennucci, ''Conformal Riemannian metrics in space of
curves'', EUSIPCO04, MIA, 2004.

\item A. Yezzi and A. Mennucci, \textquotedblright Metrics in the space of
curves\textquotedblright,\newline arXiv:math.DG/0412454, v2, May 25, 2005.
\end{enumerate}

\end{document}